\documentclass[A4,twoside,11pt,reqno]{article}
\usepackage{amssymb}
\usepackage{amsthm}
\usepackage[tbtags]{amsmath}
\usepackage{doc}
\usepackage{latexsym}
\usepackage{amscd}
\usepackage[frame,cmtip,arrow,matrix,line,graph,curve]{xy}
\usepackage{graphics}
\usepackage{epsfig}
\usepackage{eucal}
\usepackage{mathrsfs}
\newtheorem{thm}{Theorem}
\newtheorem{lem}{Lemma}
\newtheorem{exam}{Example}
\newtheorem{rem}{Remark}
\newtheorem{defn}{Definition}
\newcommand{\R}{\mathbb{R}}
\newcommand{\N}{\mathbb{N}}
\font\headd=cmr8
\allowdisplaybreaks
\begin{document}
\thispagestyle{plain}
 \markboth{}{}
\small{\addtocounter{page}{0} \pagestyle{plain}}
\begin{center}
\noindent{\large\bf Existence and Stability of Fractional Differential Equations Involving Generalized Katugampola Derivative}
\end{center}\vspace{-0.1in}
\footnote{{}\\ \\[-0.5cm]
* Corresponding Author.\\
2010 Mathematics Subject Classification: 26A33; 34K37; 35B35.\\
Key words and phrases: Fractional differential equations; Fixed point theory; Stability of solutions.\\
}
\begin{center}
\noindent{\sc Sandeep P Bhairat$^1*$ and D B Dhaigude$^{2}$}
\newline
{\it $^{1,2}$Department of Mathematics, Dr Babasaheb Ambedkar Marathwada University,\\ Aurangabad--431 004, (M.S.) India.\\
e-mail} :{\verb|sandeeppb7@gmail.com|} {\it and} {\verb|dnyanraja@gmail.com|}
\vspace{0.15in}\\
\end{center}
{\footnotesize {\sc Abstract.} The existence and stability results for a class of fractional differential equations involving generalized Katugampola derivative are presented herein. Some fixed point theorems are used and enlightening examples of obtained result are also given.}
\pagestyle{myheadings}
 \markboth{\headd Sandeep P Bhairat$~~~~~~~~~~~~~~~~~~~~~~~~~~~~~~~~~~~~~~~~~~~~~~~~~\,$}
 {\headd $~~~~~~~~~~~$ Existence and Stability of F D E Involving Generalized Katugampola Derivative}\\

\section{Introduction}
The study on various qualitative properties of solutions to numerous fractional differential and integral equations is the key topic of applied mathematics research. Indeed, every fractional differential operator is defined through the corresponding fractional integral, they are nonlocal in nature and hence more applicable than traditional one. On the other hand, fractional calculus has proven to be an useful tool in the description of various complex phenomena in the real world problems. During this theoretical development of the calculus of arbitrary order, plenty of fractional integral and differential operators are introduced and/or used by timely mathematicians, see \cite{kst} and references \cite{an,aj},\cite{db1}-\cite{db4},\cite{hh}-\cite{kh}. Wyel, Liouville, Riemann-Liouville (R-L), Hadamard are few of them. Although the well-developed theory and many more applications of the said operators, still this is a spotlight area of research in applied sciences.

Recently, U. Katugampola generalized the above mentioned fractional integral and differential operators in \cite{ki,kd}. In the same work, he obtained boundedness of generalized fractional integral in an extended Lebesgue measurable space and illustrative examples are given. In \cite{ke}, he obtained existence and uniqueness results to the solution of (initial value problem) IVP for a class of generalized fractional differential equations (FDEs).

R. Almeida, et.al \cite{am} derived these results with its Caputo counterpart and the formula for numerical treatment to solve the corresponding FDE is given. R. Almeida \cite{ra1} studied certain problems of calculus of variations dependent on Lagrange function with the same approach for first and second order. In 2015, D. Anderson et.al. \cite{du} studied the properties of the Katugampola fractional derivative with potential application in quantum mechanics. They constructed a Hamiltonian from its self adjoint operator and applied to the particle in a box model.

Recently, D.S.Oliveira, et.al.\cite{oo} proposed a generalization of Katugampola and Caputo-Katugampola fractional derivatives with the name Hilfer-Katugampola fractional derivative. This new fractional derivative interpolates the well-known fractional derivatives: Hilfer, Hilfer-Hadamard, Katugampola, Caputo-Katugampola, Riemann-Liouville, Hadamard, Caputo, Caputo-Hadamard, Liouville, Wyel as its particular cases. Following the results of \cite{kmf}, they further obtained existence and uniqueness of solution of nonlinear FDEs involving this generalized Katugampola derivative with initial condition.

The stability of functional equations was first posed by Ulam \cite{smu}. Thereafter, this type of stability evolved as an interesting field of research. The concept of stability of functional equations arises when the functional equation is being replaced by an inequality which acts as a perturbation of the functional equation, see the monograph \cite{smj} and the references cited therein. The considerable attention paid to recent development of stability results for FDEs can be found in \cite{ab,as,bl,cc,db5,kmf,smj},\cite{tm}-\cite{wlz}.

In the present work, we initiate to study the stability of the following IVPs for generalized Katugampola fractional differential equations:
\begin{equation}\label{p1}
\begin{cases}
\big({^\rho{D}_{a+}^{\alpha,\beta}x} \big)(t)&=f\big( t,x(t),({^\rho{D}_{a+}^{\alpha,\beta}x})(t)\big);\qquad t\in\Omega,\\
\big({^\rho{I}_{a+}^{1-\gamma}x}\big)(a)&=c_1,\qquad c_1\in\R,\gamma=\alpha+\beta(1-\alpha),
\end{cases}
\end{equation}
where $\alpha\in(0,1),\beta\in[0,1],\rho>0,f:\Omega\times\R\times\R\to\R$ is the given function and
\begin{equation}\label{p2}
\begin{cases}
\big({^\rho{D}_{a+}^{\alpha,\beta}x} \big)(t)&=f(t,x(t));\qquad\qquad\qquad t\in\Omega,\\
\big({^\rho{I}_{a+}^{1-\gamma}x}\big)(a)&=c_2,\qquad c_2\in\R,\gamma=\alpha+\beta(1-\alpha),
\end{cases}
\end{equation}
where $\alpha\in(0,1),\beta\in[0,1],\rho>0,f:\Omega\times\R\to\R$ is the given function, and ${^\rho{D}_{a+}^{\alpha,\beta}}, {^\rho{I}_{a+}^{1-\gamma}}$ are the generalized Katugampola fractional derivative (of order $\alpha$ and type $\beta$) and Katugampola fractional integral (of order $1-\gamma$) with $a>0,$ respectively.

The rest of the paper is organised as follows: in Section 2 we recall some preliminary facts that we need in the sequel. In Section 3 we present our main results on existence and stability of considered problems. As an application of main results, two illustrative examples are given in last Section.

\section{Preliminaries}
Let $\Omega=[a,b](0<a<b<\infty).$ As usual $C$ denotes the Banach space of all continuous functions $x:\Omega\to E$ with the superemum (uniform) norm
\begin{equation*}
{\|x\|}_{\infty}=\sup_{t\in\Omega}{\|x(t)\|}_{E}
\end{equation*}
and $AC(\Omega)$ be the space of absolutely continuous functions from $\Omega$ into $E.$ Denote $AC^1(\Omega)$ the space defined by
\begin{equation*}
AC^1(\Omega)=\bigg\{x:\Omega\to E|\frac{d}{dt}x(t)\in AC(\Omega)\bigg\}.
\end{equation*}
Let $\delta_{\rho}=t^{\rho-1}\frac{d}{dt}$ and throughout the paper, $n=[\alpha]+1,$ mention $[\alpha]$ as the integer part of $\alpha.$ Define the space
\begin{equation*}
AC_{\delta_\rho}^{n}=\big\{x:\Omega\to E|{\delta_\rho^{n-1}}x(t)\in AC(\Omega)\big\},\quad n\in\N.
\end{equation*}
Here we define the weighted space of continuous functions $g$ on ${\Omega}^*=(a,b]$ by
\begin{equation*}
C_{\gamma,\rho}(\Omega)=\bigg\{g:{\Omega}^*\to\R|\big(\frac{t^\rho-a^\rho}{\rho}\big)^{1-\gamma}g(t)\in C(\Omega)\bigg\},\quad0<\gamma\leq1,
\end{equation*}
with the norm
\begin{equation*}
{\|g\|}_{C_{\gamma,\rho}}=\big{\|{\big(\frac{t^\rho-a^\rho}{\rho}\big)}^{1-\gamma}g(t)\big\|}_{C}=\max_{t\in\Omega}{\big|{\big(\frac{t^\rho-a^\rho}{\rho}\big)}^{1-\gamma}g(t)\big|}
\end{equation*}
and
\begin{equation*}
C_{\delta_{\rho},\gamma}^{1}(\Omega)=\{g\in C(\Omega):\delta_\rho g\in C_{\gamma,\rho}(\Omega)\}
\end{equation*}
with the norms
\begin{equation*}
{\|g\|}_{C_{\delta_{\rho},\gamma}^{1}}={\|g\|}_{C}+{\|\delta_\rho g\|}_{C_{\gamma,\rho}}\quad\text{and}\quad {\|g\|}_{C_{\delta_{\rho}}^{1}}=\sum_{k=0}^{1}\max_{t\in\Omega}\big|{\delta}_\rho^k g(t)\big|.
\end{equation*}
Note that $C_{\delta_{\rho},\gamma}^{0}(\Omega)=C_{\delta_{\rho},\gamma}(\Omega)$ and $C_{0,\rho}(\Omega)=C(\Omega).$

Now we recollect definitions, some useful results and properties from  \cite{as,ki,kd,oo}:
\begin{defn}[Katugampola fractional integral]\label{ki}
Let $\alpha\in{\R}_+,c\in\R$ and $g\in{X_{c}^{p}(a,b)},$ where ${X_{c}^{p}(a,b)}$ is the space of Lebesgue measurable functions. The Katugampola fractional integral of order $\alpha$ is defined by
\begin{equation*}
({^\rho{I}_{a+}^{\alpha}g})(t)=\int_{a}^{t}s^{\rho-1}{\bigg(\frac{t^\rho-s^\rho}{\rho}\bigg)}^{\alpha-1} \frac{g(s)}{\Gamma(\alpha)}ds,\qquad t>a,\rho>0,
\end{equation*}
where $\Gamma(\cdot)$ is a Euler's gamma function.
\end{defn}
\begin{defn}[Katugampola fractional derivative]\label{kd}
Let $\alpha\in{\R}_+\setminus\N$ and $\rho>0.$ The Katugampola fractional derivative ${^\rho{D}_{a+}^{\alpha,\beta}}$ of order $\alpha$ is defined by
\begin{align*}
({^\rho{D}_{a+}^{\alpha}g})(t)&=\delta_\rho^n{({^\rho{I}_{a+}^{n-\alpha}g})(t)}\\
&={\bigg(t^{1-\rho}\frac{d}{dt}\bigg)}^n\int_{a}^{t}s^{\rho-1}{\bigg(\frac{t^\rho-s^\rho}{\rho}\bigg)}^{n-\alpha-1}\frac{g(s)}{\Gamma(n-\alpha)}ds.
\end{align*}
\end{defn}
\begin{defn}[Caputo-Katugampola fractional derivative]\label{ckd}
The Caputo-Katugampola fractional derivative $({^\rho_C{D}_{a+}^{\alpha}g})(t)$ is defined via the above Katugampola fractional derivative Definition \ref{kd} as follows
\begin{equation*}
({^\rho_C{D}_{a+}^{\alpha}g})(t)=\big({^\rho{D}_{a+}^{\alpha}}[g(s)-g(a)]\big)(t),\quad 0<\alpha<1.
\end{equation*}
\end{defn}
\begin{defn}[Generalized Katugampola fractional derivative]
The generalized Katugampola fractional derivative of order $\alpha\in(0,1)$ and type $\beta\in[0,1]$ with respect to $t$ and is defined by
\begin{align}\label{gk}
({^\rho{D}_{a\pm}^{\alpha,\beta}g})(t)=(\pm{^\rho{I}_{a\pm}^{\beta(1-\alpha)}\delta_\rho{^\rho{I}_{a\pm}^{(1-\beta)(1-\alpha)}g}})(t),\quad\rho>0
\end{align}
for the functions for which right hand side expression exists. Clearly, $0<(1-\beta)(1-\alpha)<1.$
\end{defn}
\begin{rem}\label{r1}
The generalized Katugampola operator ${^\rho{D}_{a+}^{\alpha,\beta}}$ can be written in terms of Katugampola fractional derivative as
\begin{equation*}
{^\rho{D}_{a+}^{\alpha,\beta}}={^\rho{I}_{a+}^{\beta(1-\alpha)}}{\delta_\rho}{^\rho{I}_{a+}^{1-\gamma}}={^\rho{I}_{a+}^{\beta(1-\alpha)}}{^\rho{D}_{a+}^{\gamma}},\quad \gamma=\alpha+\beta-\alpha\beta.
\end{equation*}
\end{rem}
\begin{lem}
[Semigroup property]
Let $\alpha,\beta>0,1\leq p\leq\infty,0<a<b<\infty$ and $\rho,c\in\R$ for $\rho\geq c.$ Then, for $g\in{X_{c}^{p}(a,b)}$ the following relation hold:
\begin{equation*}
({^\rho{I}_{a+}^{\alpha}}{^\rho{I}_{a+}^{\beta}g})(t)=({^\rho{I}_{a+}^{\alpha+\beta}g})(t).
\end{equation*}
\end{lem}
\begin{lem}\label{l2}
Let $t>a,~{^\rho{I}_{a+}^{\alpha}}$ and ${^\rho{D}_{a+}^{\alpha}}$ are as in Definition \ref{ki} and Definition \ref{kd}, respectively. Then the following hold:
\begin{align*}
(i)&~~\bigg({^\rho{I}_{a+}^{\alpha}}{\big(\frac{s^\rho-a^\rho}{\rho}\big)}^{\sigma}\bigg)(t)=\frac{\Gamma(\sigma+1)}{\Gamma(\sigma+\alpha+1)}{\big(\frac{t^\rho-a^\rho}{\rho}\big)}^{\sigma+\alpha},\quad\alpha\geq0,\sigma>0,\\
(ii)&~~\text{for}~\sigma=0,~~\bigg({^\rho{I}_{a+}^{\alpha}}{\big(\frac{s^\rho-a^\rho}{\rho}\big)}^{\sigma}\bigg)(t)=\big({^\rho{I}_{a+}^{\alpha}}1\big)(t)=\frac{{\big(\frac{t^\rho-a^\rho}{\rho}\big)}^{\alpha}}{\Gamma(\alpha+1)},\quad\alpha\geq0,\\
(iii)&~~\text{for}~0<\alpha<1,~~\bigg({^\rho{D}_{a+}^{\alpha}}{\big(\frac{s^\rho-a^\rho}{\rho}\big)}^{\alpha-1}\bigg)(t)=0.
\end{align*}
\end{lem}

\section{Main results}
In this section, we present the results on the existence, the attractivity and the Ulam stability of solutions for functional fractional differential equations involving generalized Katugampola fractional derivatives.

Denote $BC=BC(I),~I=[a,\infty).$ Let $D(\neq\phi)\subset BC,$ and let $G:D\to D.$ Consider the solutions of equation
\begin{equation}\label{m1}
(Gx)(t)=x(t).
\end{equation}
We define the attractivity of solutions for equation \eqref{m1} as follows:
\begin{defn}
A solutions of equation \eqref{m1} are locally attractive if there exists a ball $B(x_0,\mu)$ in the space $BC$ such that, for arbitrary solutions $y=y(t)$ and $z=z(t)$ of equation \eqref{m1} belonging to $B(x_0,\mu)\cap D,$ we have
\begin{equation}\label{m2}
\lim_{t\to\infty}(y(t)-z(t))=0.
\end{equation}
Whenever the limit \eqref{m2} is uniform with respect to $B(x_0,\mu)\cap D,$ solutions of equation \eqref{m1} are said to be uniformly locally attractive.
\end{defn}
\begin{lem}\cite{cc}\label{lm}
Let $X\subset BC.$ Then $X$ is relatively compact in $BC$ if the following conditions hold:\\
1. $X$ is uniformly bounded in $BC,$\\
2. The functions belonging to $X$ are almost equicontinuous on ${\R}_+,$ i.e. equicontinuous on every compact of ${\R}_+,$\\
3. The functions from $X$ are equiconvergent, i.e. given $\epsilon>0$ there corresponds $T(\epsilon)>0$ such that $|x(t)-\lim_{t\to\infty}x(t)|<\epsilon$ for any $t\geq T(\epsilon)$ and $x\in X.$
\end{lem}
Now we discuss the existence and the attractivity of solutions of IVP \eqref{p1}. Throughout the work, we mean $BC_{\gamma,\rho}=BC_{\gamma,\rho}(I)$ is a weighted space of all bounded and continuous functions defined by
\begin{equation*}
BC_{\gamma,\rho}=\bigg\{ x:(a,\infty]\to\R~|~{\big(\frac{t^\rho-a^\rho}{\rho}\big)}^{1-\gamma}x(t)\in BC\bigg\}
\end{equation*}
with the norm
\begin{equation*}
{\|x\|}_{BC_{\gamma,\rho}}=\sup_{t\in {\R}_+}\bigg|{\big(\frac{t^\rho-a^\rho}{\rho}\big)}^{1-\gamma}x(t)\bigg|.
\end{equation*}
\begin{thm}\cite{gd}\label{sfpt}[Schauder fixed point theorem]
Let $E$ be a Banach space and $Q$ be a nonempty bounded convex and closed subset of $E$ and $\Lambda:Q\to Q$ is compact, and continuous map. Then $\Lambda$ has at least one fixed point in $Q.$
\end{thm}
\begin{defn}
A solution of problem \eqref{p1} is a measurable function $x\in BC_{\gamma,\rho}$ satisfying initial condition $({^\rho{I}_{a+}^{1-\gamma}x})(a^+)=c_1$ and the equation $\big({^\rho{D}_{a+}^{\alpha,\beta}x} \big)(t)=f\big(t,x(t),({^\rho{D}_{a+}^{\alpha,\beta}x})(t)\big)$ on $I.$
\end{defn}
From (\cite{oo}, Theorem 3 pp. 9), we conclude the following lemma.
\begin{lem}\label{equi}
Let $\gamma=\alpha+\beta(1-\alpha),$ where $0<\alpha<1,0\leq\beta\leq1$ and $\rho>0.$ Let $f:\Omega\times\R\times\R\to\R$ be such that $f(\cdot,x(\cdot),y(\cdot))\in BC_{\gamma,\rho}$ for any $x,y\in BC_{\gamma,\rho}.$ Then problem \eqref{p1} is equivalent to Volterra integral equation of kind
\begin{equation*}
x(t)=\frac{c_1}{\Gamma(\gamma)}{\bigg(\frac{t^\rho-a^\rho}{\rho}\bigg)}^{\gamma-1}+\int_{a}^{t}s^{\rho-1}{\bigg(\frac{t^\rho-s^\rho}{\rho}\bigg)}^{\alpha-1}\frac{g(s)}{\Gamma(\alpha)}ds,
\end{equation*}
where $g(\cdot)\in{BC_{\gamma,\rho}}$ such that
\begin{equation}\label{g}
g(t)=f\bigg(t,~~\frac{c_1}{\Gamma(\gamma)}{\big(\frac{t^\rho-a^\rho}{\rho}\big)}^{\gamma-1}+({^\rho{I}_{a+}^{\alpha}g})(t),~~g(t)\bigg).
\end{equation}
\end{lem}
We use the following hypotheses in the sequel:
\begin{description}
\item[$(H_1)$] The function $t\longmapsto f(t,x,y)$ is measurable on $I$ for each $x,y\in{BC_{\gamma,\rho}},$ and the functions $x\longmapsto f(t,x,y)$ and $y\longmapsto f(t,x,y)$ are continuous on $BC_{\gamma,\rho}$ for a.e. $t\in I;$
\item[$(H_2)$] There exists a continuous function $p:I\to{\R_+}$ such that
\begin{equation*}
|f(t,x,y)|\leq\frac{p(t)}{1+|x|+|y|},~~\text{for a.e.}~t\in I,~\text{and each}~~x,y\in\R.
\end{equation*}
\end{description}
Moreover, assume that
\begin{equation*}
\lim_{t\to\infty}{\big(\frac{t^\rho-a^\rho}{\rho}\big)}^{1-\gamma}\big({^\rho{I}_{a+}^{\alpha}p}\big)(t)=0.
\end{equation*}
Set
\begin{equation*}
p^*=\sup_{t\in I}{\big(\frac{t^\rho-a^\rho}{\rho}\big)}^{1-\gamma}\big({^\rho{I}_{a+}^{\alpha}p}\big)(t).
\end{equation*}
\begin{thm}\label{t2}
Assume that the hypotheses $(H_1)$ and $(H_2)$ hold. Then the problem \eqref{p1} has at least one solution defined on $I.$ Moreover, solutions of porblem \eqref{p1} are locally attractive.
\end{thm}
\begin{proof}
For any $x\in{BC_{\gamma,\rho}},$ define the operator $\Lambda$ such that
\begin{equation}\label{p3}
(\Lambda{x})(t)=\frac{c_1}{\Gamma(\gamma)}{\bigg(\frac{t^\rho-a^\rho}{\rho}\bigg)}^{\gamma-1}+\int_{a}^{t}s^{\rho-1}{\bigg(\frac{t^\rho-s^\rho}{\rho}\bigg)}^{\alpha-1}\frac{g(s)}{\Gamma(\alpha)}ds,
\end{equation}
where $g\in{BC_{\gamma,\rho}}$ given by \eqref{g}. The operator $\Lambda$ is well defined and maps $BC_{\gamma,\rho}$ into $BC_{\gamma,\rho}.$ Indeed, the map $\Lambda(x)$ is continuous on $I$ for any $x\in{BC_{\gamma,\rho}},$ and for each $t\in{I},$ we have
\begin{align*}
\bigg|{\bigg(\frac{t^\rho-a^\rho}{\rho}\bigg)}^{1-\gamma}(\Lambda{x})(t)\bigg|
&\leq\frac{|c_1|}{\Gamma(\gamma)}+{\bigg(\frac{t^\rho-a^\rho}{\rho}\bigg)}^{1-\gamma}\int_{a}^{t}s^{\rho-1}{\bigg(\frac{t^\rho-s^\rho}{\rho}\bigg)}^{\alpha-1}\frac{|g(s)|}{\Gamma(\alpha)}ds\\
&\leq\frac{|c_1|}{\Gamma(\gamma)}+\frac{{\big(\frac{t^\rho-a^\rho}{\rho}\big)}^{1-\gamma}}{\Gamma(\alpha)}\int_{a}^{t}s^{\rho-1}{\bigg(\frac{t^\rho-s^\rho}{\rho}\bigg)}^{\alpha-1}p(s)ds\\
&\leq\frac{|c_1|}{\Gamma(\gamma)}+{\bigg(\frac{t^\rho-a^\rho}{\rho}\bigg)}^{1-\gamma}\big({^\rho{I}_{a+}^{\alpha}p}\big)(t).
\end{align*}
Thus
\begin{equation}\label{p4}
{\|\Lambda(x)\|}_{BC_{\gamma,\rho}}\leq\frac{|c_1|}{\Gamma(\gamma)}+p^*:=M.
\end{equation}
Hence, $\Lambda(x)\in{BC_{\gamma,\rho}}.$ This proves that operator $\Lambda$ maps $BC_{\gamma,\rho}$ into itself.

By Lemma \ref{equi}, the problem of finding the solutions of problem \eqref{p1} is reduced to the finding solution of the operator equation $\Lambda(x)=x.$ Equation \eqref{p4} implies that $\Lambda$ transforms the ball $B_M:=B(0,M)=\{x\in{BC_{\gamma,\rho}}:{\|x\|}_{BC_{\gamma,\rho}}\leq{M}\}$ into itself.

Now we show that the operator $\Lambda$ satisfies all the assumptions of Theorem \ref{sfpt}. The proof is given in following steps:\\
Step 1: $\Lambda$ is continuous.

Let $\{x_n\}_{n\in{\N}}$ be a sequence such that $x_n\to{x}$ in $B_M.$ Then, for each $t\in{I},$ we have
\begin{align}\label{p5}
\bigg|{\bigg(\frac{t^\rho-a^\rho}{\rho}\bigg)}^{1-\gamma}&(\Lambda{x_n})(t)-{\bigg(\frac{t^\rho-a^\rho}{\rho}\bigg)}^{1-\gamma}(\Lambda{x})(t)\bigg|\nonumber\\
&\leq\frac{{\big(\frac{t^\rho-a^\rho}{\rho}\big)}^{1-\gamma}}{\Gamma(\alpha)}\int_{a}^{t}s^{\rho-1}{\bigg(\frac{t^\rho-s^\rho}{\rho}\bigg)}^{\alpha-1}|g_n(s)-g(s)|ds,
\end{align}
where $g,g_n\in{BC_{\gamma,\rho}},$ $g$ is defined by \eqref{g} and
\begin{equation*}
g_n(t)=f\bigg(t,~\frac{c_1}{\Gamma(\gamma)}{\bigg(\frac{t^\rho-a^\rho}{\rho}\bigg)}^{\gamma-1}+({^\rho{I}_{a+}^{\alpha}g_n})(t),~g_n(t)\bigg).
\end{equation*}
If $t\in{I},$ then from \eqref{p5}, we obtain
\begin{align}\label{p6}
\bigg|{\bigg(\frac{t^\rho-a^\rho}{\rho}\bigg)}^{1-\gamma}&(\Lambda{x_n})(t)-{\bigg(\frac{t^\rho-a^\rho}{\rho}\bigg)}^{1-\gamma}(\Lambda{x})(t)\bigg|\nonumber\\
&\leq2\frac{{\big(\frac{t^\rho-a^\rho}{\rho}\big)}^{1-\gamma}}{\Gamma(\alpha)}\int_{a}^{t}s^{\rho-1}{\bigg(\frac{t^\rho-s^\rho}{\rho}\bigg)}^{\alpha-1}p(s)ds,
\end{align}
Since $x_n\to x$ as $n\to\infty$ and ${\big(\frac{t^\rho-a^\rho}{\rho}\big)}^{1-\gamma}\big({^\rho{I}_{a+}^{\alpha}p}\big)(t)\to0$ as $t\to\infty$ then \eqref{p6} implies
\begin{equation*}
{\|\Lambda(x_n)-\Lambda(x)\|}_{BC_{\gamma,\rho}}\to0\quad\text{as }n\to\infty.
\end{equation*}
Step 2: $\Lambda(B_M)$ is uniformly bounded.

This is clear since $\Lambda(B_M)\subset{B_M}$ and ${B_M}$ is bounded.\\
Step 3: $\Lambda(B_M)$ is equicontinuous on every compact subset $[a,T]$ of $I,~T>a.$

Let $t_1,t_2\in[a,T],t_1<t_2$ and $x\in B_M.$ We have
\begin{align*}
\bigg|{\bigg(\frac{{t_2}^\rho-a^\rho}{\rho}\bigg)}^{1-\gamma}&(\Lambda{x})(t_2)-{\bigg(\frac{{t_1}^\rho-a^\rho}{\rho}\bigg)}^{1-\gamma}(\Lambda{x})(t_1)\bigg|\nonumber\\
&\leq\bigg|\frac{{\big(\frac{{t_2}^\rho-a^\rho}{\rho}\big)}^{1-\gamma}}{\Gamma(\alpha)}\int_{a}^{{t_2}}s^{\rho-1}{\bigg(\frac{{t_2}^\rho-s^\rho}{\rho}\bigg)}^{\alpha-1}g(s)ds\\
&~~-\frac{{\big(\frac{{t_1}^\rho-a^\rho}{\rho}\big)}^{1-\gamma}}{\Gamma(\alpha)}\int_{a}^{{t_1}}s^{\rho-1}{\bigg(\frac{{t_1}^\rho-s^\rho}{\rho}\bigg)}^{\alpha-1}g(s)ds\bigg|,
\end{align*}
with $g(\cdot)\in{BC_{\gamma,\rho}}$ given by \eqref{g}. Thus we get
\begin{align*}
\bigg|{\bigg(\frac{{t_2}^\rho-a^\rho}{\rho}\bigg)}^{1-\gamma}&(\Lambda{x})(t_2)-{\bigg(\frac{{t_1}^\rho-a^\rho}{\rho}\bigg)}^{1-\gamma}(\Lambda{x})(t_1)\bigg|\nonumber\\
&\leq\frac{{\big(\frac{{t_2}^\rho-a^\rho}{\rho}\big)}^{1-\gamma}}{\Gamma(\alpha)}\int_{t_1}^{{t_2}}s^{\rho-1}{\bigg(\frac{{t_2}^\rho-s^\rho}{\rho}\bigg)}^{\alpha-1}|g(s)|ds\\
&~~+\int_{a}^{t_1}\bigg|\bigg[{\bigg(\frac{{t_2}^\rho-a^\rho}{\rho}\bigg)}^{1-\gamma}s^{\rho-1}{\bigg(\frac{{t_2}^\rho-s^\rho}{\rho}\bigg)}^{\alpha-1}\\
&~~-{\bigg(\frac{{t_1}^\rho-a^\rho}{\rho}\bigg)}^{1-\gamma}s^{\rho-1}{\bigg(\frac{{t_1}^\rho-s^\rho}{\rho}\bigg)}^{\alpha-1}\bigg]\bigg|\frac{|g(s)|}{\Gamma(\alpha)}ds\\
&\leq\frac{{\big(\frac{{t_2}^\rho-a^\rho}{\rho}\big)}^{1-\gamma}}{\Gamma(\alpha)}\int_{t_1}^{{t_2}}s^{\rho-1}{\bigg(\frac{{t_2}^\rho-s^\rho}{\rho}\bigg)}^{\alpha-1}p(s)ds\\
&~~+\int_{a}^{t_1}\bigg|\bigg[{\bigg(\frac{{t_2}^\rho-a^\rho}{\rho}\bigg)}^{1-\gamma}s^{\rho-1}{\bigg(\frac{{t_2}^\rho-s^\rho}{\rho}\bigg)}^{\alpha-1}\\
&~~-{\bigg(\frac{{t_1}^\rho-a^\rho}{\rho}\bigg)}^{1-\gamma}s^{\rho-1}{\bigg(\frac{{t_1}^\rho-s^\rho}{\rho}\bigg)}^{\alpha-1}\bigg]\bigg|\frac{p(s)}{\Gamma(\alpha)}ds.
\end{align*}
Thus, for $p_*=\sup_{t\in[a,T]}p(t)$ and from the continuity of the function $p,$ we obtain
\begin{align*}
\bigg|{\bigg(\frac{{t_2}^\rho-a^\rho}{\rho}\bigg)}^{1-\gamma}&(\Lambda{x})(t_2)-{\bigg(\frac{{t_1}^\rho-a^\rho}{\rho}\bigg)}^{1-\gamma}(\Lambda{x})(t_1)\bigg|\nonumber\\
&\leq{p_*}\frac{{\big(\frac{{T}^\rho-a^\rho}{\rho}\big)}^{1-\gamma+\alpha}}{\Gamma(\alpha+1)}{\bigg(\frac{t_2^\rho-t_1^\rho}{\rho}\bigg)}^{\alpha}\\
&~~+\frac{p_*}{\Gamma(\alpha)}\int_{a}^{t_1}\bigg|{\bigg(\frac{{t_2}^\rho-a^\rho}{\rho}\bigg)}^{1-\gamma}s^{\rho-1}{\bigg(\frac{{t_2}^\rho-s^\rho}{\rho}\bigg)}^{\alpha-1}\\
&~~-{\bigg(\frac{{t_1}^\rho-a^\rho}{\rho}\bigg)}^{1-\gamma}s^{\rho-1}{\bigg(\frac{{t_1}^\rho-s^\rho}{\rho}\bigg)}^{\alpha-1}\bigg|ds.
\end{align*}
As $t_1\to{t_2},$ the right hand side of the above inequation tends to zero.\\
Step 4: $\Lambda(B_M)$ is equiconvergent.

Let $t\in{I}$ and $x\in{B_M},$ then we have
\begin{equation*}
\bigg|{\bigg(\frac{{t}^\rho-a^\rho}{\rho}\bigg)}^{1-\gamma}(\Lambda{x})(t)\bigg|\leq\frac{|c_1|}{\Gamma(\gamma)}+\frac{{\big(\frac{t^\rho-a^\rho}{\rho}\big)}^{1-\gamma}}{\Gamma(\alpha)}\int_{a}^{t}s^{\rho-1}{\bigg(\frac{t^\rho-s^\rho}{\rho}\bigg)}^{\alpha-1}|g(s)|ds
\end{equation*}
where $g(\cdot)\in{BC_{\gamma,\rho}}$ is given by \eqref{g}. Thus we get
\begin{align*}
\bigg|{\bigg(\frac{{t}^\rho-a^\rho}{\rho}\bigg)}^{1-\gamma}(\Lambda{x})(t)\bigg|&\leq\frac{|c_1|}{\Gamma(\gamma)}+\frac{{\big(\frac{t^\rho-a^\rho}{\rho}\big)}^{1-\gamma}}{\Gamma(\alpha)}\int_{a}^{t}s^{\rho-1}{\bigg(\frac{t^\rho-s^\rho}{\rho}\bigg)}^{\alpha-1}p(s)ds\\
&\leq\frac{|c_1|}{\Gamma(\gamma)}+{\bigg(\frac{t^\rho-a^\rho}{\rho}\bigg)}^{1-\gamma}\big({^\rho{I}_{a+}^{\alpha}p}\big)(t).
\end{align*}
Since ${\big(\frac{t^\rho-a^\rho}{\rho}\big)}^{1-\gamma}\big({^\rho{I}_{a+}^{\alpha}p}\big)(t)\to0$ as $t\to\infty,$ then, we get
\begin{equation*}
|(\Lambda{x})(t)|\leq\frac{|c_1|}{{\big(\frac{t^\rho-a^\rho}{\rho}\big)}^{1-\gamma}\Gamma(\gamma)}+\frac{{\big(\frac{t^\rho-a^\rho}{\rho}\big)}^{1-\gamma}\big({^\rho{I}_{a+}^{\alpha}p}\big)(t)}{{\big(\frac{t^\rho-a^\rho}{\rho}\big)}^{1-\gamma}}\to0~~\text{as}~~t\to\infty.
\end{equation*}
Hence
\begin{equation*}
|(\Lambda{x})(t)-(\Lambda{x})(+\infty)|\to0~~\text{as }t\to\infty.
\end{equation*}
In view of Lemma \ref{lm} and immediate consequence of Steps 1 to 4, we conclude that $\Lambda:B_M\to{B_M}$ is continuous and compact. Theorem \ref{sfpt} implies that $\Lambda$ has a fixed point $x$ which is a solution of problem \eqref{p1} on $I.$\\
Step 5: Local attactivity of solutions.

Let $x_0$ is a solution of IVP \eqref{p1}. Taking $x\in{B(x_0,2p^*)},$ we have
\begin{align*}
\bigg|{\bigg(\frac{{t}^\rho-a^\rho}{\rho}\bigg)}^{1-\gamma}(\Lambda{x})(t)&-{\bigg(\frac{{t}^\rho-a^\rho}{\rho}\bigg)}^{1-\gamma}{x_0}(t)\bigg|\\
&=\bigg|{\bigg(\frac{{t}^\rho-a^\rho}{\rho}\bigg)}^{1-\gamma}(\Lambda{x})(t)-{\bigg(\frac{{t}^\rho-a^\rho}{\rho}\bigg)}^{1-\gamma}(\Lambda{x_0})(t)\bigg|\\
&\leq\frac{{\big(\frac{{t}^\rho-a^\rho}{\rho}\big)}^{1-\gamma}}{\Gamma(\alpha)}\int_{a}^{t}s^\rho{\bigg(\frac{t^\rho-s^\rho}{\rho}\bigg)}^{\alpha-1}|f(s,g(s))-f(s,g_0(s))|ds,
\end{align*}
where $g,g_0\in{BC_{\gamma,\rho}},~g$ is given by \eqref{g} and
\begin{equation*}
g_0(t)=f\bigg(t,~\frac{c_1}{\Gamma(\gamma)}{\bigg(\frac{t^\rho-a^\rho}{\rho}\bigg)}^{\gamma-1}+({^\rho{I}_{a+}^{\alpha}g_0})(t),~g_0(t)\bigg).
\end{equation*}
Then
\begin{align*}
\bigg|{\bigg(\frac{{t}^\rho-a^\rho}{\rho}\bigg)}^{1-\gamma}(\Lambda{x})(t)-{\bigg(\frac{{t}^\rho-a^\rho}{\rho}\bigg)}^{1-\gamma}{x_0}(t)\bigg|&\leq2\frac{{\big(\frac{{t}^\rho-a^\rho}{\rho}\big)}^{1-\gamma}}{\Gamma(\alpha)}\int_{a}^{t}s^\rho{\bigg(\frac{t^\rho-s^\rho}{\rho}\bigg)}^{\alpha-1}p(s)ds\\
&\leq2p^*.
\end{align*}
We obtain
\begin{equation*}
{\|(\Lambda({x})-x_0\|}_{BC_{\gamma,\rho}}\leq2p^*.
\end{equation*}
Hence $\Lambda$ is a continuous function such that $\Lambda(B(x_0,2p^*))\subset{B(x_0,2p^*)}.$

Moreover, if $x$ is a solution of IVP \eqref{p1}, then
\begin{align*}
|x(t)-x_0(t)|&=|(\Lambda{x})(t)-(\Lambda{x_0})(t)|\\
&\leq\frac{1}{\Gamma(\alpha)}\int_{a}^{t}s^\rho{\bigg(\frac{t^\rho-a^\rho}{\rho}\bigg)}^{\alpha-1}|g(s)-g_0(s)|ds\\
&\leq2({^\rho{I}_{a+}^{\alpha}p})(t).
\end{align*}
Thus
\begin{equation}\label{p7}
|x(t)-x_0(t)|\leq2\frac{{\big(\frac{{t}^\rho-a^\rho}{\rho}\big)}^{1-\gamma}({^\rho{I}_{a+}^{\alpha}p})(t)}{{\big(\frac{{t}^\rho-a^\rho}{\rho}\big)}^{1-\gamma}}.
\end{equation}
With the fact that $\lim_{t\to\infty}{{\big(\frac{{t}^\rho-a^\rho}{\rho}\big)}^{1-\gamma}}({^\rho{I}_{a+}^{\alpha}p})(t)=0$ and inequation \eqref{p7}, we obtain
\begin{equation*}
\lim_{t\to\infty}|x(t)-x_0(t)|=0.
\end{equation*}
Consequently, all solutions of \eqref{p1} are locally attractive.
\end{proof}
Now onwards in this section, we deal with the existence and the Ulam stability of solutions for problem \eqref{p2}.
\begin{lem}
Let $\gamma=\alpha+\beta(1-\alpha),$ where $0<\alpha<1,0\leq\beta\leq1$ and $\rho>0.$ Let $f:\Omega\times\R\to\R$ be such that $f(\cdot,x(\cdot))\in{C_{\gamma,\rho}(\Omega)}$ for any $x\in{C_{\gamma,\rho}(\Omega)}.$ Then problem \eqref{p2} is equivalent to the Volterra integral equation of kind
\begin{equation*}
x(t)=\frac{c_2}{\Gamma(\gamma)}{{\big(\frac{{t}^\rho-a^\rho}{\rho}\big)}}^{\gamma-1}+({^\rho{I}_{a+}^{\alpha}}f(\cdot,x(\cdot)))(t).
\end{equation*}
\end{lem}
Let $\epsilon>0$ and $\Phi:\Omega\to[0,\infty)$ be a continuous function and consider the following inequalities:
\begin{equation}\label{a1}
|\big({^\rho{D}_{a+}^{\alpha,\beta}x}\big)(t)-f(t,x(t))|\leq\epsilon; \qquad t\in\Omega,
\end{equation}
\begin{equation}\label{a2}
|\big({^\rho{D}_{a+}^{\alpha,\beta}x}\big)(t)-f(t,x(t))|\leq\Phi(t); \qquad t\in\Omega,
\end{equation}
\begin{equation}\label{a3}
|\big({^\rho{D}_{a+}^{\alpha,\beta}x}\big)(t)-f(t,x(t))|\leq\epsilon\Phi(t); \qquad t\in\Omega.
\end{equation}
\begin{defn}\label{d7}
Problem \eqref{p2} is Ulam-Hyers stable if there exists a real number $\psi>0$ such that for each $\epsilon>0$ and for each solution $x\in{C_{\gamma,\rho}}$ of inequality \eqref{a1} there exists a solution $\bar{x}\in{C_{\gamma,\rho}}$ of problem \eqref{p2} with
\begin{equation*}
|x(t)-\bar{x}(t)|\leq\epsilon{\psi};\qquad t\in\Omega.
\end{equation*}
\end{defn}
\begin{defn}\label{d8}
Problem \eqref{p2} is generalized Ulam-Hyers stable if there exists $\Psi:C([0,\infty),[0,\infty))$ with $\Psi(0)=0$ such that for each $\epsilon>0$ and for each solution $x\in{C_{\gamma,\rho}}$ of inequality \eqref{a1} there exists a solution $\bar{x}\in{C_{\gamma,\rho}}$ of problem \eqref{p2} with
\begin{equation*}
|x(t)-\bar{x}(t)|\leq \Psi(\epsilon);\qquad t\in\Omega.
\end{equation*}
\end{defn}
\begin{defn}\label{d9}
Problem \eqref{p2} is Ulam-Hyers-Rassias stable with respect to $\Phi$ if there exists a real number $\psi_\phi>0$ such that for each $\epsilon>0$ and for each solution $x\in{C_{\gamma,\rho}}$ of inequality \eqref{a3} there exists a $\bar{x}\in{C_{\gamma,\rho}}$ of problem \eqref{p2} with
\begin{equation*}
|x(t)-\bar{x}(t)|\leq\epsilon{\psi_\phi}\Phi(t);\qquad t\in\Omega.
\end{equation*}
\end{defn}
\begin{defn}\label{d10}
Problem \eqref{p2} is generalized Ulam-Hyers-Rassias stable with respect to $\Phi$ if there exists a real number $\psi_\phi>0$ such that for each solution $x\in{C_{\gamma,\rho}}$ of inequality \eqref{a2} there exists a $\bar{x}\in{C_{\gamma,\rho}}$ of problem \eqref{p2} with
\begin{equation*}
|x(t)-\bar{x}(t)|\leq{\psi_\phi}\Phi(t);\qquad t\in\Omega.
\end{equation*}
\end{defn}
\begin{rem}
It is clear that
\begin{description}
  \item[(i)] Definition \ref{d7} $\Rightarrow$ Definition \ref{d8}.
  \item[(ii)] Definition \ref{d9} $\Rightarrow$ Definition \ref{d10}.
  \item[(iii)] Definition \ref{d9} for $\Phi(\cdot)=1~ \Rightarrow$ Definition \ref{d8}.
\end{description}
\end{rem}
\begin{defn}
A solution of problem \eqref{p2} is a measurable function $x\in{C_{\gamma,\rho}}$ that satisfies the condition $\big({^\rho{I}_{a+}^{1-\gamma}x}\big)(a)=c_2,$ and the differential equation $\big({^\rho{D}_{a+}^{\alpha,\beta}x} \big)(t)=f(t,x(t))$ on $\Omega.$
\end{defn}
Now we introduce the following hypotheses which will be used in the sequel:
\begin{description}
  \item[$(H_3)$] The function $t\longmapsto{f(t,x)}$ is measurable on $\Omega$ for each $x\in{C_{\gamma,\rho}},$ and the function $x\longmapsto{f(t,x)}$ is continuous on ${C_{\gamma,\rho}}$ for a.e. $t\in\Omega,$
  \item[$(H_4)$] There exists a continuous function $p:\Omega\to[0,\infty)$ such that
\begin{equation*}
|f(t,x)|\leq\frac{p(t)}{1+|x|}|x|,~~\text{for a.e. }t\in\Omega,~~\text{and each }x\in\R.
\end{equation*}
\end{description}
Set $p^*=\sup_{t\in\Omega}p(t).$ Now we shall give the existence theorem in the following:
\begin{thm}\label{t3}
Assume that the hypotheses $(H_3)$ and $(H_4)$ hold. Then the problem \eqref{p2} has at least one solution defined on $\Omega.$
\end{thm}
\begin{proof}
Consider the operator $\Lambda:C_{\gamma,\rho}\to{C_{\gamma,\rho}}$ such that
\begin{equation}\label{p8}
(\Lambda{x})(t)=\frac{c_2}{\Gamma(\gamma)}{\bigg(\frac{t^\rho-a^\rho}{\rho}\bigg)}^{\gamma-1}+\int_{a}^{t}s^{\rho-1}{\bigg(\frac{t^\rho-s^\rho}{\rho}\bigg)}^{\alpha-1}\frac{f(s,x(s))}{\Gamma(\alpha)}ds.
\end{equation}
Clearly, the fixed points of this operator equation $(\Lambda{x})(t)=x(t)$ are solutions of IVP \eqref{p2}. For any $x\in{C_{\gamma,\rho}}$ and each $t\in\Omega,$ we have
\begin{align*}
\bigg|{\bigg(\frac{t^\rho-a^\rho}{\rho}\bigg)}^{1-\gamma}(\Lambda{x})(t)\bigg|
&\leq\frac{|c_2|}{\Gamma(\gamma)}+{\bigg(\frac{t^\rho-a^\rho}{\rho}\bigg)}^{1-\gamma}\int_{a}^{t}s^{\rho-1}{\bigg(\frac{t^\rho-s^\rho}{\rho}\bigg)}^{\alpha-1}\frac{|f(s,x(s))|}{\Gamma(\alpha)}ds\\
&\leq\frac{|c_2|}{\Gamma(\gamma)}+\frac{{\big(\frac{t^\rho-a^\rho}{\rho}\big)}^{1-\gamma}}{\Gamma(\alpha)}\int_{a}^{t}s^{\rho-1}{\bigg(\frac{t^\rho-s^\rho}{\rho}\bigg)}^{\alpha-1}p(s)ds\\
&\leq\frac{|c_2|}{\Gamma(\gamma)}+\frac{p^*}{\Gamma(\alpha)}{\big(\frac{t^\rho-a^\rho}{\rho}\big)}^{1-\gamma}\int_{a}^{t}s^{\rho-1}{\bigg(\frac{t^\rho-s^\rho}{\rho}\bigg)}^{\alpha-1}ds\\
&\leq\frac{|c_2|}{\Gamma(\gamma)}+\frac{p^*}{\Gamma(\alpha+1)}{\big(\frac{t^\rho-a^\rho}{\rho}\big)}^{1-\gamma}{\big(\frac{t^\rho-a^\rho}{\rho}\big)}^{\alpha}\\
&\leq\frac{|c_2|}{\Gamma(\gamma)}+\frac{p^*}{\Gamma(\alpha+1)}{\big(\frac{T^\rho-a^\rho}{\rho}\big)}^{\alpha+1-\gamma}.
\end{align*}
Thus
\begin{equation}\label{p9}
{\|\Lambda{x}\|}_{C}\leq\frac{|c_2|}{\Gamma(\gamma)}+\frac{p^*}{\Gamma(\alpha+1)}{\big(\frac{T^\rho-a^\rho}{\rho}\big)}^{\alpha+1-\gamma}:=N.
\end{equation}
Thus $\Lambda$ transforms the ball $B_{N}=B(0,N)=\{z\in{C_{\gamma,\rho}}:{\|z\|}_{C}\leq{N}\}$ into itself. We shall show that the operator $\Lambda:B_N\to{B_N}$ satisfies all the conditions of Theorem \ref{abc}. The proof is given in following several steps.\\
Step 1: $\Lambda:B_N\to{B_N}$ is continuous.

Let $\{x_n\}_{n\in{\N}}$ be a sequence such that $x_n\to{x}$ in $B_N.$ Then, for each $t\in{I},$ we have
\begin{align}\label{p10}
\bigg|{\bigg(\frac{t^\rho-a^\rho}{\rho}\bigg)}^{1-\gamma}&(\Lambda{x_n})(t)-{\bigg(\frac{t^\rho-a^\rho}{\rho}\bigg)}^{1-\gamma}(\Lambda{x})(t)\bigg|\nonumber\\
&\leq\frac{{\big(\frac{t^\rho-a^\rho}{\rho}\big)}^{1-\gamma}}{\Gamma(\alpha)}\int_{a}^{t}s^{\rho-1}{\bigg(\frac{t^\rho-s^\rho}{\rho}\bigg)}^{\alpha-1}|f(s,x_n(s))-f(s,x(s))|ds.
\end{align}
Since $x_n\to{x}$ as $n\to\infty$ and $f$ is continuous, then by Lebesgue dominated convergence theorem, inequation \eqref{p10} implies ${\|\Lambda(x_n)-\Lambda(x)\|}_{C}\to0$ as $n\to\infty.$\\
Step 2: $\Lambda(B_N)$ is uniformly bounded.

Since $\Lambda(B_N)\subset{B_N}$ and $B_N$ is bounded. Hence, $\Lambda(B_N)$ is uniformly bounded.\\
Step 3: $\Lambda(B_N)$ is equicontinuous.

Let $t_1,t_2\in\Omega,t_1<t_2$ and $x\in B_N.$ We have
\begin{align*}
\bigg|{\bigg(\frac{{t_2}^\rho-a^\rho}{\rho}\bigg)}^{1-\gamma}&(\Lambda{x})(t_2)-{\bigg(\frac{{t_1}^\rho-a^\rho}{\rho}\bigg)}^{1-\gamma}(\Lambda{x})(t_1)\bigg|\\
&\leq\bigg|\frac{{\big(\frac{{t_2}^\rho-a^\rho}{\rho}\big)}^{1-\gamma}}{\Gamma(\alpha)}\int_{a}^{{t_2}}s^{\rho-1}{\bigg(\frac{{t_2}^\rho-s^\rho}{\rho}\bigg)}^{\alpha-1}f(s,x(s))ds\\
&~~-\frac{{\big(\frac{{t_1}^\rho-a^\rho}{\rho}\big)}^{1-\gamma}}{\Gamma(\alpha)}\int_{a}^{{t_1}}s^{\rho-1}{\bigg(\frac{{t_1}^\rho-s^\rho}{\rho}\bigg)}^{\alpha-1}f(s,x(s))ds\bigg|\\
&\leq\frac{{\big(\frac{{t_2}^\rho-a^\rho}{\rho}\big)}^{1-\gamma}}{\Gamma(\alpha)}\int_{t_1}^{{t_2}}s^{\rho-1}{\bigg(\frac{{t_2}^\rho-s^\rho}{\rho}\bigg)}^{\alpha-1}|f(s,x(s))|ds\\
&~~+\int_{a}^{t_1}\bigg|\bigg[{\bigg(\frac{{t_2}^\rho-a^\rho}{\rho}\bigg)}^{1-\gamma}s^{\rho-1}{\bigg(\frac{{t_2}^\rho-s^\rho}{\rho}\bigg)}^{\alpha-1}\\
&~~-{\bigg(\frac{{t_1}^\rho-a^\rho}{\rho}\bigg)}^{1-\gamma}s^{\rho-1}{\bigg(\frac{{t_1}^\rho-s^\rho}{\rho}\bigg)}^{\alpha-1}\bigg]\bigg|\frac{|f(s,x(s))|}{\Gamma(\alpha)}ds\\
&\leq\frac{{\big(\frac{{t_2}^\rho-a^\rho}{\rho}\big)}^{1-\gamma}}{\Gamma(\alpha)}\int_{t_1}^{{t_2}}s^{\rho-1}{\bigg(\frac{{t_2}^\rho-s^\rho}{\rho}\bigg)}^{\alpha-1}p(s)ds\\
&~~+\int_{a}^{t_1}\bigg|\bigg[{\bigg(\frac{{t_2}^\rho-a^\rho}{\rho}\bigg)}^{1-\gamma}s^{\rho-1}{\bigg(\frac{{t_2}^\rho-s^\rho}{\rho}\bigg)}^{\alpha-1}\\
&~~-{\bigg(\frac{{t_1}^\rho-a^\rho}{\rho}\bigg)}^{1-\gamma}s^{\rho-1}{\bigg(\frac{{t_1}^\rho-s^\rho}{\rho}\bigg)}^{\alpha-1}\bigg]\bigg|\frac{p(s)}{\Gamma(\alpha)}ds.
\end{align*}
Thus, for $p_{*}=\sup_{t\in\Omega}p(t)$ and from the continuity of the function $p,$ we obtain
\begin{align*}
\bigg|{\bigg(\frac{{t_2}^\rho-a^\rho}{\rho}\bigg)}^{1-\gamma}&(\Lambda{x})(t_2)-{\bigg(\frac{{t_1}^\rho-a^\rho}{\rho}\bigg)}^{1-\gamma}(\Lambda{x})(t_1)\bigg|\\
&\leq\frac{p_*}{\Gamma(\alpha+1)}{{\bigg(\frac{{T}^\rho-a^\rho}{\rho}\bigg)}^{1-\gamma+\alpha}}{\bigg(\frac{t_2^\rho-t_1^\rho}{\rho}\bigg)}^{\alpha}\\
&~~+\frac{p_*}{\Gamma(\alpha)}\int_{a}^{t_1}\bigg|{\bigg(\frac{{t_2}^\rho-a^\rho}{\rho}\bigg)}^{1-\gamma}s^{\rho-1}{\bigg(\frac{{t_2}^\rho-s^\rho}{\rho}\bigg)}^{\alpha-1}\\
&~~-{\bigg(\frac{{t_1}^\rho-a^\rho}{\rho}\bigg)}^{1-\gamma}s^{\rho-1}{\bigg(\frac{{t_1}^\rho-s^\rho}{\rho}\bigg)}^{\alpha-1}\bigg|ds.
\end{align*}
As $t_1\to{t_2},$ the right hand side of the above inequality tends to zero.

As a consequence of Steps 1 to 3 together with Arzela-Ascoli Theorem, we can conclude that $\Lambda$ is continuous and compact. By applying the Schauder fixed point theorem, we conclude that $\Lambda$ has a fixed point $x$ which is a solution of the problem \eqref{p2}.
\end{proof}
\begin{thm}\label{t4}
Assume that $(H_3),~(H_4)$ and the following hypotheses hold:
\begin{description}
  \item[$(H_5)$] There exists $\lambda_\phi>0$ such that for each $t\in\Omega,$ we have
  \begin{equation*}
    ({^\rho{I}_{a+}^{\alpha}}\Phi(t))\leq\lambda_\phi\Phi(t);
  \end{equation*}
  \item[$(H_6)$] There exists $q\in{C(\Omega,[0,\infty))}$ such that for each $t\in\Omega$,
  \begin{equation*}
    p(t)\leq{q(t)}\Phi(t).
  \end{equation*}
\end{description}
Then the problem \eqref{p2} is generalized Ulam-Hyers-Rassias stable.
\end{thm}
\begin{proof}
Consider the operator $\Lambda:{C_{\gamma,\rho}}\to{C_{\gamma,\rho}}$ defined in \eqref{p8}. Let $x$ be a solution of inequality \eqref{a2}, and let us assume that $\bar{x}$ is a solution of problem \eqref{p2}. Thus
\begin{equation*}
\bar{x}(t)=\frac{c_2}{\Gamma(\gamma)}{\bigg(\frac{t^\rho-a^\rho}{\rho}\bigg)}^{\gamma-1}+\int_{a}^{t}s^{\rho-1}{\bigg(\frac{t^\rho-s^\rho}{\rho}\bigg)}^{\alpha-1}\frac{f(s,\bar{x}(s))}{\Gamma(\alpha)}ds.
\end{equation*}
From the inequality \eqref{a2} for each $t\in\Omega,$ we have
\begin{equation*}
\bigg|x(t)-\frac{c_2}{\Gamma(\gamma)}{\bigg(\frac{t^\rho-a^\rho}{\rho}\bigg)}^{\gamma-1}-\int_{a}^{t}s^{\rho-1}{\bigg(\frac{t^\rho-s^\rho}{\rho}\bigg)}^{\alpha-1}\frac{f(s,x(s))}{\Gamma(\alpha)}ds\bigg|\leq\Phi(t).
\end{equation*}
Set $q^*=\sup_{t\in\Omega}q(t).$ From the hypotheses (H5) and (H6), for each $t\in\Omega,$ we get
\begin{align*}
\big|x(t)-\bar{x}(t)\big|&\leq\bigg|x(t)-\frac{c_2}{\Gamma(\gamma)}{\bigg(\frac{t^\rho-a^\rho}{\rho}\bigg)}^{\gamma-1}-\int_{a}^{t}s^{\rho-1}{\bigg(\frac{t^\rho-s^\rho}{\rho}\bigg)}^{\alpha-1}\frac{f(s,x(s))}{\Gamma(\alpha)}ds\bigg|\\
&~~+\int_{a}^{t}s^{\rho-1}{\bigg(\frac{t^\rho-s^\rho}{\rho}\bigg)}^{\alpha-1}\frac{|f(s,x(s))-f(s,\bar{x}(s))|}{\Gamma(\alpha)}ds\\
&\leq\Phi(t)+\int_{a}^{t}s^{\rho-1}{\bigg(\frac{t^\rho-s^\rho}{\rho}\bigg)}^{\alpha-1}\frac{2q^*\Phi(s)}{\Gamma(\alpha)}ds\\
&\leq\Phi(t)+2q^*({^\rho{I}_{a+}^{\alpha}}\Phi)(t)\\
&\leq\Phi(t)+2q^*\lambda_\phi\Phi(t)\\
&=[1+2q^*\lambda_\phi]\Phi(t).
\end{align*}
Thus
\begin{equation*}
|x(t)-\bar{x}(t)|\leq{\psi_\phi\Phi(t)}.
\end{equation*}
Hence the problem \eqref{p2} is generalized Ulam-Hyers-Rassias stable.
\end{proof}
\noindent Define the metric
\begin{equation*}
d(x,y)=\sup_{t\in\Omega}\frac{{\big(\frac{t^\rho-a^\rho}{\rho}\big)}^{1-\gamma}{|x(t)-y(t)|}}{\Phi(t)}
\end{equation*}
in the space $C_{\gamma,\rho}(\Omega).$ The following fixed point theorem is used in our further result.
\begin{thm}\cite{dm}\label{abc}
Let $\Theta:{C_{\gamma,\rho}}\to{C_{\gamma,\rho}}$ be a strictly contractive operator with a Lipschitz constant $L<1.$ There exists a nonnegative integer $k$ such that $d(\Theta^{k+1}x,\Theta^kx)<\infty$ for some $x\in{C_{\gamma,\rho}},$ then the following propositions hold true:\\
(A1) The sequence ${\{\Theta^k{x}\}}_{n\in\N}$ converges to a fixed point $x^*$ of $\Theta;$\\
(A2) $x^*$ is a unique fixed point of $\Theta$ in $X=\{y\in{C_{\gamma,\rho}(\Omega)}:d(\Theta^kx,y)<\infty\};$\\
(A3) If $y\in{X},$ then $d(y,x^*)\leq\frac{1}{1-L}d(y,\Theta{x}).$
\end{thm}
\begin{thm}\label{t6}
Assume that $(H_5)$ and the following hypothesis hold:
\begin{description}
\item[$(H_7)$] There exists $\phi\in{C(\Omega,[0,\infty))}$ such that for each $t\in\Omega,$ and all $x,\bar{x}\in\R,$ we have
\begin{equation*}
|f(t,x)-f(t,\bar{x}|\leq{\big(\frac{t^\rho-a^\rho}{\rho}\big)}^{1-\gamma}\phi(t)\Phi(t)|x-\bar{x}|.
\end{equation*}
\end{description}
If
\begin{equation*}
L={\big(\frac{T^\rho-a^\rho}{\rho}\big)}^{1-\gamma}\phi^*\lambda_\phi<1,
\end{equation*}
where $\phi^*=\sup_{t\in\Omega}\phi(t),$ then there exists a unique solution $x_0$ of problem \eqref{p2}, and the problem \eqref{p2} is generalized Ulam-Hyers-Rassias stable. Furthermore, we have
\begin{equation*}
|x(t)-\bar{x}(t)|\leq\frac{\Phi(t)}{1-L}.
\end{equation*}
\end{thm}
\begin{proof}
Let $\Lambda:{C_{\gamma,\rho}}\to{C_{\gamma,\rho}}$ be the operator defined in \eqref{p8}. Apply Theorem \ref{abc}, we have
\begin{align*}
|(\Lambda{x})(t)-(\Lambda{\bar{x}})(t)|&\leq\int_{a}^{t}s^{\rho-1}{\bigg(\frac{t^\rho-s^\rho}{\rho}\bigg)}^{\alpha-1}\frac{|f(s,x(s))-f(s,\bar{x}(s))|}{\Gamma(\alpha)}ds\\
&\leq\int_{a}^{t}s^{\rho-1}{\bigg(\frac{t^\rho-s^\rho}{\rho}\bigg)}^{\alpha-1}\phi(s)\Phi(s)\frac{|{\big(\frac{s^\rho-a^\rho}{\rho}\big)}^{1-\gamma}x(s)-{\big(\frac{s^\rho-a^\rho}{\rho}\big)}^{1-\gamma}\bar{x}(s)|}{\Gamma(\alpha)}ds\\
&\leq\int_{a}^{t}s^{\rho-1}{\bigg(\frac{t^\rho-s^\rho}{\rho}\bigg)}^{\alpha-1}\phi^*(s)\Phi(s)\frac{{\|x-\bar{x}\|}_C}{\Gamma(\alpha)}ds\\
&\leq\phi^*({^\rho{I}_{a+}^{\alpha}})\Phi(t){{\|x-\bar{x}\|}_C}\\
&\leq\phi^*\lambda_\phi\Phi(t){{\|x-\bar{x}\|}_C}.
\end{align*}
Thus
\begin{equation*}
\bigg|{\bigg(\frac{t^\rho-a^\rho}{\rho}\bigg)}^{1-\gamma}(\Lambda{x})(t)-{\bigg(\frac{t^\rho-a^\rho}{\rho}\bigg)}^{1-\gamma}(\Lambda{\bar{x}})(t)\bigg|\leq{\bigg(\frac{T^\rho-a^\rho}{\rho}\bigg)}^{1-\gamma}\phi^*\lambda_\phi\Phi(t){{\|x-\bar{x}\|}_C}.
\end{equation*}
Hence
\begin{equation*}
d(\Lambda(x),\Lambda(\bar{x}))=\sup_{t\in\Omega}\frac{{\|(\Lambda{x})(t)-(\Lambda{\bar{x}})(t)\|}_C}{\Phi(t)}\leq{L{\|x-\bar{x}\|}_C}
\end{equation*}
from which we conclude the theorem.
\end{proof}
\section{Examples}
In this section we present some examples to illustrate our results.
\begin{exam}
Consider the following IVP with generalized Katugampola fractional derivative:
\begin{equation}\label{e1}
\begin{cases}
\big({^\rho{D}_{a+}^{\frac{1}{2}}x}\big)(t)&=f(t,x,y);\qquad t\in[a,b],\\
\big({^\rho{I}_{a+}^{\frac{1}{4}}x}\big)(a)&=(1-a),
\end{cases}
\end{equation}
where $\alpha=\frac{1}{2},\beta=\frac{1}{2},\rho>0,\gamma=\frac{3}{4},0<a<b\leq{e},$ and
\begin{equation*}
\begin{cases}
f(t,x,y)&=\frac{\theta{(t-a)}^{-\frac{1}{4}}\sin{(t-a)}}{64(1+\sqrt{t-a})(1+|x|+|y|)};\quad t\in(a,b],~x,y\in\R,\\
f(a,x,y)&=0;\qquad x,y\in\R.
\end{cases}
\end{equation*}
\end{exam}
\noindent Clearly, the function $f$ is continuous for each $x,y\in\R.$ The hypothesis $(H_2)$ is satisfied with
\begin{equation*}
\begin{cases}
p(t)&=\frac{\theta{(t-a)}^{-\frac{1}{4}}|\sin{(t-a)}|}{64(1+\sqrt{t-a})};\quad0<\theta\leq1,~t\in(a,+\infty),\\
p(a)&=0.
\end{cases}
\end{equation*}
Thus, all the conditions of Theorem \ref{t2} are satisfied. Hence, the problem \eqref{e1} has at least one solution defined on $[a,+\infty).$

Also, we have
\begin{align*}
{\bigg(\frac{t^\rho-a^\rho}{\rho}\bigg)}^{1-\gamma}({^\rho{I}_{a+}^{\frac{1}{2}}p})(t)&={\bigg(\frac{t^\rho-a^\rho}{\rho}\bigg)}^{\frac{1}{4}}
\int_{a}^{t}s^{\rho-1}{\bigg(\frac{t^\rho-s^\rho}{\rho}\bigg)}^{-\frac{1}{2}}\frac{p(s)}{\Gamma(\frac{1}{2})}ds\\
&\leq\frac{1}{8}{\bigg(\frac{t^\rho-a^\rho}{\rho}\bigg)}^{-\frac{1}{4}}\to0\quad\text{as}~~t\to+\infty.
\end{align*}
This implies that the solutions of problem \eqref{e1} are locally asymptotically stable.
\begin{exam}
Consider the following problem of FDE involving generalized Katugampola derivative:
\begin{equation}\label{e2}
\begin{cases}
\big({^\rho{D}_{a+}^{\frac{1}{2}}x}\big)(t)&=f(t,x);\qquad t\in[a,b],\\
\big({^\rho{I}_{a+}^{\frac{1}{4}}x}\big)(a)&=(1-a),
\end{cases}
\end{equation}
where $\alpha=\frac{1}{2},\beta=\frac{1}{2},\rho>0,\gamma=\frac{3}{4},0<a<b\leq{e},$ and
\begin{equation*}
\begin{cases}
f(t,x)&=\frac{\theta{(t-a)}^{-\frac{1}{4}}\sin{(t-a)}}{64(1+\sqrt{t-a})(1+|x|)};\quad t\in(a,b],~x\in\R,\\
f(a,x)&=0;\qquad x\in\R.
\end{cases}
\end{equation*}
\end{exam}
Clearly, the function $f$ is continuous for all $x\in\R.$ Hypothesis $(H_4)$ is satisfied with
\begin{equation*}
\begin{cases}
p(t)&=\frac{\theta{(t-a)}^{-\frac{1}{4}}|\sin{(t-a)}|}{64(1+\sqrt{t-a})};\quad0<\theta\leq1,~t\in(a,b],~x\in\R,\\
p(a)&=0.
\end{cases}
\end{equation*}
Hence, Theorem \ref{t3} implies that problem \eqref{e2} has at least one solution defined on $[a,b].$ Also hypothesis $(H_5)$ is satisfied with
\begin{equation*}
\Phi(t)=e^3, \quad\text{and}\quad \lambda_\phi=\frac{1}{\Gamma(\frac{3}{2})}.
\end{equation*}
Consequently, Theorem \ref{t4} implies that problem \eqref{e2} is generalized Ulam-Hyers-Rassias stable.

\end{document}